\documentclass[11pt,letterpaper]{amsart}
\usepackage[parfill]{parskip}
\usepackage{amssymb}
\usepackage{amsmath}
\usepackage{amsthm}
\usepackage{amsfonts}
\usepackage{graphicx}

\usepackage[T1]{fontenc}

\usepackage{textcomp}
\usepackage{palatino}
\usepackage[text={6.5in, 8.5in}, centering]{geometry}

\usepackage[usenames]{color}

\newtheorem{theorem}{Theorem}

\theoremstyle{definition}

\theoremstyle{remark}
\newtheorem{remark}[theorem]{Remark}

\def\ocirc#1{\ifmmode\setbox0=\hbox{$#1$}\dimen0=\ht0
    \advance\dimen0 by1pt\rlap{\hbox to\wd0{\hss\raise\dimen0
    \hbox{\hskip.2em$\scriptscriptstyle\circ$}\hss}}#1\else
    {\accent"17 #1}\fi}

\newcommand{\R}{\mathbb{R}}

\newcommand{\curl}{\operatorname{curl}}

\title[Vorticity in $L \log L$]{A spatially localized  $L \log L$ estimate on the vorticity in the 3D NSE}
\author{Z. Bradshaw}
\author{Z. Gruji\'c}

\begin{document}

\maketitle

\begin{abstract}
The purpose of this note is to present a spatially localized $L \log L$ bound on the
vorticity in the 3D Navier-Stokes equations, assuming a very mild, \emph{purely geometric}
condition. This yields an extra-log decay of the distribution function of the vorticity, which
in turn implies \emph{breaking the criticality} in a physically, numerically, and mathematical
analysis-motivated criticality scenario based on vortex stretching and anisotropic diffusion.
\end{abstract}

\section{Introduction}

\noindent Three-dimensional Navier-Stokes equations (3D NSE) -- describing a flow of 3D
incompressible viscous fluid -- read

\[
 u_t+(u\cdot \nabla)u=-\nabla p + \triangle u,
\]

\noindent supplemented with the incompressibility condition $ \, \mbox{div} \,
u = 0$, where $u$ is the velocity of the fluid, and $p$ is the
pressure (here, the viscosity is set to 1). Applying the curl operator yields the vorticity 
formulation,

\begin{equation}\label{vorticity}
 \omega_t+(u\cdot \nabla)\omega= (\omega \cdot \nabla)u + \triangle
 \omega,
\end{equation}

\noindent where $\omega = \, \mbox{curl} \, u$ is the vorticity of the fluid. Taking the spatial
domain to be the whole space, the velocity
can be recovered from the vorticity via the Biot-Savart Law,

\[
 u(x) = c \int \nabla\frac{1}{|x-y|} \times \omega(y) \, dy;
\]

\noindent this makes (\ref{vorticity}) a closed (non-local) system for the vorticity field alone.

An \emph{a priori} $L^1$-bound on the evolution of the vorticity in the 3D NSE was obtained 
by Constantin in \cite{Co90} . This is a vorticity analogue of 
Leray's \emph{a priori} $L^2$-bound on the velocity, and both quantities scale in the same
fashion. So far, there has been no \emph{a priori} estimate on the weak solutions to the 3D NSE breaking
this scaling.

The goal of this short article is to show that a very mild, \emph{purely geometric} assumption yields a
uniform-in-time $L \log L$ bound on the vorticity, effectively breaking the aforementioned
scaling. More precisely, the assumption is a
uniform-in-time boundedness of the localized \emph{vorticity direction} in the weighted $bmo$-space
$\displaystyle{\widetilde{bmo}_{\frac{1}{|\log r|}}}$. An interesting feature of this space (cf. \cite{Jan76}) 
is that it allows for
discontinuous functions exhibiting singularities of, e.g.,  $\sin \log |\log ( \, \mbox{something algebraic} \, )|$-type.
In particular, the vorticity direction can blow-up in a geometrically spectacular fashion
-- every point on the unit sphere being a limit point -- and the $L \log L$ bound will still hold.

Besides being of an independent interest, the $L \log L$ bound on the vorticity implies an
extra-log decay of the distribution function.
This is significant as it transforms a recently exposed \cite{Gr12, DaGr12-1}
large-data criticality scenario for the 3D NSE into a no blow-up scenario. Shortly, creation
and persistence (in the sense of the time-average) of the \emph{axial} filamentary scales comparable
to the \emph{macro-scale}, paired with the $L^1$-induced decay of the volume of the suitably defined
region of intense vorticity, leads to 
creation and persistence of the \emph{transversal micro-scales} comparable to the scale of 
local, \emph{anisotropic linear
sparseness}, enabling the anisotropic diffusion to equalize the nonlinear effects. The extra log-decay
of the volume transforms the equalizing scenario into the anisotropic diffusion-win scenario.

The present result is, in a way, complementary to the results obtained by the authors in \cite{BrGr13}. 
The class of conditions leading to an $L \log L$-bound presented in \cite{BrGr13} can be characterized
as `wild in time' with a uniform spatial (e.g., algebraic) structure, while the condition presented
here can be characterized as `wild in space' and uniform in time. As in \cite{BrGr13}, the proof is based
on an adaptation of the method exposed in \cite{Co90}, the novel component being utilization of
\emph{analytic cancelations}
in the vortex-stretching term via the Hardy space-version of the Div-Curl Lemma \cite{CoLiMeSe}, the local
version \cite{Goldberg} of $\mathcal{H}^1-BMO$ duality \cite{Fe-72, FeSt-72}, 
and the intimate connection between the $BMO$-norm
and the logarithm. While the argument in \cite{BrGr13} relied on the structure of the evolution of the scalar 
components and the result (cf. \cite{RS87,  BeijingLectures}) stating that the $BMO$-norm of the logarithm
of a polynomial is bounded \emph{independently of the coefficients}, the present argument
relies on sharp pointwise multiplier theorem in $\widetilde{bmo}$ \cite{Jan76, NaYa85, MS85} and 
Coifman-Rochberg's 
$BMO$-estimate on the logarithm of the \emph{maximal function} of a locally integrable function
\cite{CoRoch80}, the estimate being \emph{fully independent of the function} and depending only 
on the dimension of the space.

\section{An excursion to harmonic analysis}

In this section, we compile several results from harmonic analysis that will prove useful.

\medskip

\texttt{Hardy spaces} $\mathcal{H}^1$ and $\frak{h}^1$

\noindent The maximal function of a
distribution $f$ is defined as,
\[M_hf(x)=\sup_{t>0}|f*h_t(x)|, \ x\in \R^n,
\] 
where $h$ is a fixed, normalized ($\int h ~dx=1$) test function supported in the unit ball, and 
$h_t$ denotes $t^{-n}h(\cdot/t)$.

The distribution $f$ is in the Hardy space $\mathcal {H}^1$ if $\|f\|_{\mathcal {H}^1}=\|M_hf\|_1<\infty$.

The local maximal function is defined as,
\[m_hf(x)=\sup_{0<t<1}|f*h_t(x)|, \ x\in \R^n,
\]
and the distribution $f$ is in the local Hardy space $\frak h^1$ if $\|f\|_{\frak {h}^1}=\|m_hf\|_1<\infty$.

\newpage

\texttt{Div-Curl Lemma} \ (Coifman, Lions, Meyer, Semmes \cite{CoLiMeSe})

\noindent Suppose that
$E$ and $B$ are $L^2$-vector fields satisfying $\mbox{div} \,  E = \curl B
= 0$  (in the sense of distribution).  Then,
\[
\|E\cdot B\|_{\mathcal H^1} \leq c(n) \, \|E\|_{L^2} \|B\|_{L^2}.
\]

\medskip

$BMO$ \texttt{and weighted} $bmo$ \texttt{spaces}

The classical space of bounded mean oscillations, $BMO$ is defined as follows,
\[
 BMO = \biggl\{ f \in L^1_{loc} : \, \sup_{x \in \mathbb{R}^3, r>0}
 \Omega \bigl(f, I(x,r)\bigr) < \infty \biggr\},
\]
where $\displaystyle{\Omega \bigl(f, I(x,r)\bigr)=\frac{1}{|I(x,r)|}\int_{I(x,r)} |f(x)-f_I | \, dx}$
is the mean oscillation of the function $f$ with respect to its mean
$f_I = \frac{1}{|I(x,r)|}\int_{I(x,r)} f(x) \, dx$, over the cube $I(x,r)$ centered at $x$ with
the side-length $r$.

A local version of $BMO$, usually denoted by $bmo$, is defined by finiteness of the following
expression,
\[
 \|f\|_{bmo} = \sup_{x \in \mathbb{R}^3, 0 < r < \delta} \Omega \bigl(f, I(x,r)\bigr)
 + \sup_{x \in \mathbb{R}^3, r \ge \delta} \frac{1}{|I(x,r)|} \int_{I(x,r)} |f(y)| \, dy,
\]
for some positive $\delta$.

If $f \in L^1$, we can focus on small scales, e.g., $0 < r < \frac{1}{2}$.
Let $\phi$ be a positive, non-decreasing function on $(0, \frac{1}{2})$,
and consider the following version of local
weighted spaces of bounded mean oscillations,
\[
 \|f\|_{\widetilde{bmo}_\phi} = \|f\|_{L^1} +  \sup_{x \in \mathbb{R}^3, 0<r<\frac{1}{2}}
 \frac{\Omega \bigl(f, I(x,r)\bigr)}{\phi(r)}
\]
(cf. \cite{MS85}).

Of special interest will be the spaces $\widetilde{bmo} = \widetilde{bmo}_1$, and 
$\widetilde{bmo}_{\frac{1}{|\log r|}}$.

\medskip

$\mathcal{H}^1-BMO$ \texttt{and} $\frak{h}^1-bmo$ \, \texttt{duality} \ (Fefferman \cite{Fe-72, FeSt-72}, 
Goldberg \cite{Goldberg})

\noindent $\bigl(\mathcal{H}^1\bigr)^* = BMO$ and $\bigl(\frak{h}^1\bigr)^* = bmo$; 
the duality is realized via integration of one object against the other.

\medskip

\texttt{Pointwise multipliers in $\widetilde{bmo}$} \, (\cite{Jan76, NaYa85, MS85})

\emph{A sharp pointwise multiplier theorem.} (\cite{MS85}) \ Let $h$ be in $\widetilde{bmo}$, and $g$ in 
$L^\infty \cap \widetilde{bmo}_\frac{1}{|\log r|}$. Then,
\[
 \|g \, h\|_{\widetilde{bmo}} \le c(n) \, \Bigl( \|g\|_\infty + \|g\|_{\widetilde{bmo}_\frac{1}{|\log r|}} \Bigr) 
 \ \|h\|_{\widetilde{bmo}}.
\]
More precisely, the space of pointwise $\widetilde{bmo}$ multipliers coincides with $L^\infty \cap
\widetilde{bmo}_\frac{1}{|\log r|}$.

\newpage

\texttt{Coifman and Rochberg's estimate on} $\| \log M f\|_{BMO}$  (\cite{CoRoch80})

Let $M$ denote
the Hardy-Littlewood maximal operator. Coifman and
Rochberg \cite{CoRoch80} obtained a characterization of $BMO$ in terms of images of the 
logarithm of the maximal function of non-negative locally integrable functions (plus
a bounded part). The main ingredient in demonstrating one direction is the following estimate,
\[
 \| \log M f \|_{BMO} \le c(n),
\]
for any locally integrable function $f$. (The bound is completely
\emph{independent of} $f$.)

\medskip

This estimate remains valid if we replace $M f$ with $\mathcal{M} f = \bigl( M \sqrt{|f|}\bigr)^2$ (cf. \cite{IwVe-99});
the advantage of working with $\mathcal{M}$ is that the $L^2$-maximal theorem implies
the following estimate,
\begin{equation}\label{meep}
 \|\mathcal{M} f\|_1 \le c(n) \|f\|_1,
\end{equation}
a bound that does not hold for the original maximal operator $M$.

\section{Setting and the result}

Consider a weak (distributional) Leray solution $u$ on $\mathbb{R}^3$. The vorticity 
analogue of the Leray's \emph{a priori} bound on the energy was presented
in \cite{Co90}: assuming that the initial vorticity is in $L^1$ (or, more generally, a bounded measure),
the $L^1$-norm of the vorticity remains bounded on any finite time-interval.

Our goal is to obtain a spatially localized $L \log L$ bound on the vorticity under a suitable
assumption on the structural blow-up of the \emph{vorticity direction} 
$\displaystyle{\xi=\frac{\omega}{|\omega|}}$.

Fix a spatial ball $B(0,R_0)$, and consider a test function $\psi$ 
supported in $B=B(0, 2R_0)$ such that $\psi = 1$ on $B(0, R_0)$, and $\displaystyle{
|\nabla\psi(x)| \le c \frac{1}{R_0} \psi^\delta(x)}$ for some $\delta > 0$.

Let $w=\sqrt{1+|\omega|^2}$. We aim to control the evolution of $\psi \, w \log w$; by
the Stein's lemma \cite{St-L-log-L}, this is essentially equivalent to controlling the $L^1$-norm
of $M w$.

For simplicity of the exposition, we assume that the initial vorticity is also
in $L^2$, and that $T>0$ is the first (possible) blow-up time. This way, the solution
in view is smooth on $(0,T)$, and we can focus on obtaining a $\sup_{t \in (0,T)}$-bound.
Alternatively, one can employ the retarded mollifiers.

\begin{theorem}
Let $u$ be a Leray solution to the 3D NSE. Assume that the initial vorticity $\omega_0$ 
is in $L^1 \cap L^2$, and that $T>0$ is the first (possible) blow-up time. Suppose that
\[
 \sup_{t \in (0,T)} \| (\psi \xi) (\cdot, t)  \|_{\widetilde{bmo}_\frac{1}{|\log r|}} < \infty.
\]
Then,
\[
 \sup_{t \in (0,T)} \int \psi(x) \, w(x,t) \log w(x,t) \, dx < \infty.
\]
\end{theorem}

\begin{remark}
Since $\omega_0$ is in $L^1$, in addition to the Leray's \emph{a priori} bounds on $u$,
\[
 \sup_t \|u(\cdot, t)\|_{L^2} < \infty \ \ \mbox{and} \ \ \int_t \int_x |\nabla u(x,t)|^2 \, dx \, dt < \infty,
\]
the following \emph{a priori} bounds on $\omega$ are also at our disposal \cite{Co90},
\[
 \sup_t \|\omega(\cdot, t)\|_{L^1} < \infty \ \ \mbox{and} \ \ \int_t \int_x |\nabla 
 \omega(x,t)|^\frac{4}{3+\epsilon} \, dx \, dt < \infty.
\]
\end{remark}

\begin{proof} \ Setting $q(y)=\sqrt{1+|y|^2}$, the evolution of $w=\sqrt{1+|\omega|^2}$ 
satisfies the following partial differential inequality (\cite{Co90}),

\begin{equation}\label{w}
\partial_t w - \triangle w + (u \cdot \nabla) w \le \omega \cdot \nabla u \cdot 
\frac{\omega}{w}.
\end{equation}

Since our goal is to control the evolution of $\psi \, w \log w$, it will prove convenient to
multiply (\ref{w}) by $\psi \, (1+\log w)$. Here is the calculus corresponding to each
of the four terms.

\bigskip

\texttt{time-derivative}

\[
 \partial_t w \times \psi \, (1+\log w) = \partial_t (\psi \, w \log w).
\]

\medskip

\texttt{Laplacian}

\begin{align*}
 -\triangle w \, \times & \, \psi \, (1+\log w)\\
  &= -\triangle (\psi \, w \log w) + \triangle \psi \, w \log w\\
  &+ \psi \frac{1}{w} \sum_i (\partial_i w)^2 + 2 \sum_i \partial_i \psi \, \partial_i w \, (1+\log w).
\end{align*}
  
\medskip

\texttt{advection}

\begin{align*}
(u \cdot \nabla) w \, \times & \, \psi \, (1+\log w)\\
  &= \sum_i u_i \, \partial_i w \, \psi \, (1+\log w)\\
  &=\sum_i \bigl( \partial_i (u_i \, w \, \psi \, (1+\log w)) - 
   u_i \, w \, \partial_i \psi \, (1+\log w) - u_i \, \psi \, \partial_i w\bigr)\\
  &=\sum_i \bigl( \partial_i (u_i \, w \, \psi \, (1+\log w)) - 
   u_i \, w \, \partial_i \psi \, (1+\log w) - \partial_i (u_i \psi w)
   + (u_i \, \partial_i \psi \, w)\bigr).
\end{align*}

\medskip

\texttt{vortex-stretching}

\[
 \omega \cdot \nabla u \cdot \frac{\omega}{w} \, \times \, \psi \, (1+\log w)
 = \omega \cdot \nabla u \cdot \psi \, \frac{\omega}{|\omega|}  \, (1+\log w)
 + \omega \cdot \nabla u \cdot \psi \, \biggl(\frac{\omega}{w}-\frac{\omega}{|\omega|} \biggr) 
 \, (1+\log w).
\]

\medskip

Integrating over the space-time, the above representation yields (dropping the zero and the
positive terms, and estimating the remaining terms in the straightforward fashion via
H\"older and Sobolev),

\[
I(\tau) \equiv \int \psi(x) \, w(x,\tau) \log w(x,\tau) \, dx \le 
 I(0) + c \int_0^\tau \int_x \omega \cdot \nabla u \cdot \psi \, \xi  \, \log w \, dx \, dt
 + \ \mbox{\emph{a priori} \, bounded},
\]
for any $\tau$ in $[0,T)$.

In order to take the advantage of the Coifman-Rochberg's estimate, we decompose
the logarithmic factor as 

\[
\log w = \log \frac{w}{\mathcal{M}w} 
+ \log \mathcal{M}w.
\]

Denoting $\displaystyle{\int_0^\tau \int_x \omega \cdot 
\nabla u \cdot \psi \, \xi  \, \log w \, dx \, dt}$
by $J$, this yields $J=J_1+J_2$ where

\[
 J_1 =  \int_0^\tau \int_x \omega \cdot \nabla u \cdot \psi \, \xi  \,
 \log \frac{w}{\mathcal{M}w}  \, dx \, dt
\]

and

\[
 J_2 =  \int_0^\tau \int_x \omega \cdot \nabla u \cdot \psi \, \xi \,
  \log \mathcal{M}w \, dx \, dt.
\]

For $J_1$, we use the \emph{pointwise} inequality

\[
 w \log \frac{w}{\mathcal{M}w} \le \mathcal{M}w - w
\]

(a consequence of the pointwise inequality $\mathcal{M}f \ge f$, and the inequality
$e^{x-1} \ge x$ for $x \ge 1$).

This leads to 

\[
 J_1 \le \int_0^\tau \int_x |\nabla u| \Bigl(\mathcal{M}w - w\Bigr)  \psi \, dx \, dt
\]

which is \emph{a priori} bounded by the Cauchy-Schwarz and the $L^2$-maximal theorem.

For $J_2$, we have the following string of inequalities,

\[
\begin{aligned}
J_2
 & \le c \int_0^\tau \|\omega \cdot \nabla u\|_{\frak{h}^1} 
 \|\psi \, \xi \log \mathcal{M}w\|_{bmo} \, dt\\
 & \le c \int_0^\tau \|\omega \cdot \nabla u\|_{\mathcal{H}^1} 
 \|\psi \, \xi \log \mathcal{M}w\|_{\widetilde{bmo}} \, dt\\
  & \le c \int_0^\tau \|\omega\|_2 \|\nabla u\|_2 \Bigl( \|\psi \, \xi\|_\infty + \|\psi \, \xi\|_{\widetilde{bmo}_{\frac{1}{|\log r|}}}
  \Bigr)
 \Bigl(\|\log \mathcal{M}w\|_{BMO} + \|\log \mathcal{M}w\|_1\Bigr) \, dt\\
 & \le c \sup_{t \in (0,T)} \ \biggl\{ \Bigl(1 + \|\psi \, \xi\|_{\widetilde{bmo}_{\frac{1}{|\log r|}}}\Bigr)  
  \Bigl(\|\log \mathcal{M}w\|_{BMO} + \|\log \mathcal{M}w\|_1\Bigr)
     \biggr\} \ \ \int_t \int_x |\nabla u|^2 \, dx \, dt\\ 
 &\le c \ \Bigl(1+ \sup_{t \in (0,T)} \|\omega\|_1\Bigr) \ 
 \Bigl(1 + \sup_{t \in (0,T)} \|\psi \, \xi\|_{\widetilde{bmo}_{\frac{1}{|\log r|}}}\Bigr)   
 \ \int_t \int_x |\nabla u|^2 \, dx \, dt\\
\end{aligned}
\]

by $\frak{h}^1 - bmo$ duality, the Div-Curl Lemma, the pointwise $\widetilde{bmo}$-multiplier 
theorem, the Coifman-Rochberg's estimate, and the bound (\ref{meep}) combined with a
couple of elementary inequalities. 
This completes the proof of the $L \log L$
estimate. 

\end{proof}

\subsection*{Acknowledgements.}

Z.B. acknowledges the support of the \emph{Virginia Space Grant
Consortium} via the Graduate Research Fellowship; Z.G. acknowledges
the support of the \emph{Research Council of Norway} via the grant
213474/F20, and the \emph{National Science Foundation} via the grant
DMS 1212023.
We are grateful to Luong Dang Ky for the comments on the previous version of
the paper.

\bibliographystyle{plain}
\bibliography{ref}

\end{document}